\documentclass[a4paper,12pt]{amsart}
\usepackage[latin1]{inputenc}      
\usepackage{amssymb,epsfig}

\newtheorem{theorem}{Theorem}[section]

\newtheorem{corollary}[theorem]{Corollary}

\renewenvironment{proof}{\noindent {\bf Proof.}}{ \hfill\qed\\ }

\begin{document}

\def\lrarrow{\hbox to 30pt{\ \rightarrowfill\ }}
\def\rfn#1{\smash{\mathop{\lrarrow}\limits^{#1}}}
\def\diagram#1{\vbox{\halign
     {$\hfil \displaystyle{##}\hfil $&&$\hfil \displaystyle{##}\hfil $\cr
                               #1\crcr}}}
\def\llrarrow{\hbox to 30pt{\kern1pt\rightarrowfill\kern1pt}}
\def\mrfn#1{\smash{\mathop{\llrarrow}\limits^{#1}}}
\def\lllarrow{\hbox to 30pt{\kern1pt\leftarrowfill\kern1pt}}
\def\mlfn#1{\smash{\mathop{\lllarrow}\limits^{#1}}}
\def\ufn#1{\vbox{\vskip4pt
                 \hbox{{\Big\uparrow
                        \rlap{$\vcenter{\hbox{$\scriptstyle{#1}$}}$}}}
                 \vskip-1pt}}
\def\dfn#1{\vbox{\vskip4pt
                 \hbox{{\Big\downarrow
                        \rlap{$\vcenter{\hbox{$\scriptstyle{#1}$}}$}}}
                 \vskip-1pt}}
\def\pull{\vphantom{\Big\downarrow}\raise.5em\hbox{\kern6pt\bf pull}}

\def\d{\delta}
\def\F{\mathcal{F}}
\def\U{\mathcal{U}}
\def\E{\mathcal{E}}
\def\m{\mu}
\def\e{\varepsilon}
\def\t{\theta}
\def\Z{\mathbb{Z}}
\def\R{\mathbb{R}}
\def\N{\mathbb{N}}
\def\B{\mathcal{B}}
\def\D{\mathcal{D}}
\def\M{\mathcal{M}}
\def\Bh{\hat{\mathcal{B}}}
\def\X{\hat{X}}
\def\mh{\hat{\mu}}
\def\Th{\hat{T}}
\def\hA{\hat{A}}
\def\A{\mathcal{A}}
\def\C{\mathcal{C}}
\def\Gd{\mathcal{G}_{\delta}}
\def\T2{T^{\times 2}}
\def\2m{\mu^{\times 2}}
\def\S{{\bf S^1}}

\def\sS{\mathcal{H}}
\def\e{\epsilon}
\def\w{\omega}
\def\a{\alpha}
\def\Am{\mathcal{A}_{\mu,\w}}
\def\Bm{\mathcal{B}_{\mu,t}}
\def\t{\theta}
\def\Xt{X_{\t}}
\def\Dt{D_{\t}}\def\Rt{R_{\t}}\def\Ut{U_{\t}}
\def\Bt{F_{\t}}
\def\Cs{\mathcal{C}_s}
\def\Es{\mathcal{E}_s}
\def\Ds{\mathcal{D}_s}

\title{Approximation and  billiards}
\author{Serge Troubetzkoy}

\address{Centre de physique théorique, Institut de mathématiques de Luminy Federation de Recherches des Unites de Mathematique de Marseille and
Université de la Méditerranée, Luminy, Case 907, F-13288 Marseille Cedex 9, France.}
\email{troubetz@iml.univ-mrs.fr}
\urladdr{http://iml.univ-mrs.fr/{\lower.7ex\hbox{\~{}}}troubetz/}

\begin{abstract}
This survey is based on a series of talks I gave at 
the conference ``Dynamical systems and Diophantine approximation''
at l'Instut Henri Poincaré in June 2003.  
I will present asymptotic results (transitivity, ergodicity,
weak-mixing) for billiards based on the
approximation technique developed by Katok and Zemlyakov.  
I will also present approximation techniques which allow
to prove the abundance of periodic trajectories in certain
irrational polygons.
\end{abstract}

\maketitle

\section{Introduction.}

The use of approximation techniques in dynamical systems was apparently
started in 1941 by Oxtoby and Ulam \cite{OU} who proved that  for a
finite-dimensional compact manifold with a non-atomic measure which is
positive on open sets the set of ergodic measure-preserving
homeomorphisms is generic in the strong topology.  In 1967
Katok and Stepin \cite{KS} proved the genericity of
of ergodicity and weak-mixing for certain classes of interval 
exchange transformations.
The approximation method 
was first applied to polygonal billiard by 
Katok and Zemlyakov who, using the fact that for rational polygons
the directional billiard is minimal in almost all directions, 
proved that the typical irrational
billiard is topological transitive (i.e.~has dense orbits) \cite{KZ}.

In this survey we will first reprove the Katok Zemlyakov result and then
give further developments based on this idea.  Then
we will discuss another approximation result, based on
inhomogeneous diophantine approximation, to conclude that there
are many periodic billiard orbits in irrational 
parallelograms and related polygons.

\section{Background.}

In this section we give the necessary background on billiards,
more details can be found in the surveys \cite{T},\cite{MT}.
The {\em billiard flow} in  a domain  $P \subset \R^2$ is defined as
follows:
a point mass moves freely inside $P$ and when
it reaches the boundary it is reflected following the usual law
of geometric optics, the angle of incidence equals the angle of reflection.
The first return map to the boundary of $P$ is called the {\em billiard
map}. 
The billiard map preserves a natural measure $\sin \theta \, ds \, d\theta$
where $s$ is the arc length parameter on $\partial P$ and $\theta \in
\S$ is the angle between the direction of the billiard motion and
the direction of the boundary.

A polygon is called {\em rational} is all the angles between sides
are rational multiples of $\pi$.  For rational polygons there 
is a well known construction of {\em invariant billiard surfaces}.
Suppose the angles are of the form $m_i \pi/n_i$.  
Consider the group $G(P)$ generated by the reflections in the sides
of $P$ and let $O(P) \subset \S$ be the linearization of
$G(P)$.  The group $O(P)$ is simply the dihedral group $D_N$
generated by the 
reflections in lines through the origin that meet at angles $\pi/N$.
Here $N$ is the least common multiple of the $n_i's$.

Consider the {\em phase space} of the billiard flow $P \times \S$, and let
$R_{\theta}$ be its subset of points whose second coordinate belongs
to the orbit of $\theta$ under $D_N$.  Since every trajectory
changes directions at reflections by elements of $D_N$ we have
that $R_\theta$ is an invariant set.  In fact it is a flat surface
(possibly) with  conical singularities.

To construct $R_\theta$ (for $\theta \ne k\pi/N$) consider $2N$
disjoint parallel copies $P_1,P_2,\dots,P_{2N}$ of $P$.  Consider
the $O(P)$ orbits $\{\theta_1,\theta_2,\dots,\theta_{2N}\}$ of $\theta$.
Consider a pair $(P_i,\theta_i)$.  Fix a side of $P_i$ and reflect  
$\theta_i$ in this side, this yields a $\theta_j$.  Paste the given
side of $P_i$ and $P_j$.  Doing this for each $i$ yields the surface.
For example the billiard in the square has a torus as invariant
surface while the invariant surface for the billiard in the in the
right triangle with angle $\pi/8$ is of genus two (see Figure 1).

\begin{figure}
\centerline{\psfig{file=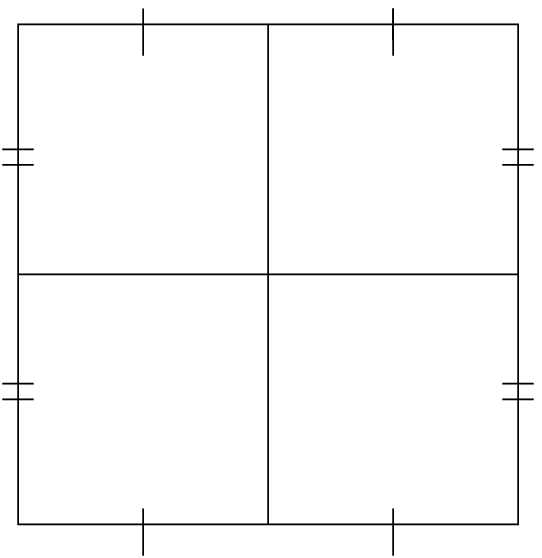,height=45mm} \hfill \psfig{file=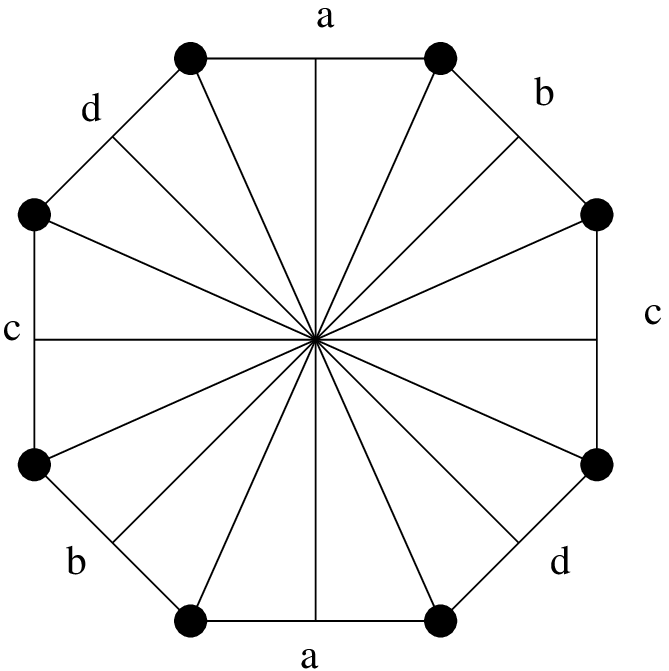,height=50mm}}
\caption{Invariant surfaces.}
\end{figure}  

\section{Old approximation results for polygons.}

It is not very difficult to prove that the billiard flow/map 
restricted to $R_{\theta}$ is minimal for all but countably many $\theta$.
We are now ready to state and prove Katok-Zemlyakov's result.
Consider the space $X$ of simply connected $n$-gons.  Note that
this space is not compact since $n$-gons can degenerate into $n'$ gons
with $n' < n$.  Thus we say that a property is {\em typical} if for
any compact set $K \subset X$ the property holds for a dense
$G_{\delta}$ subset of $K$.{\footnote{This detail is often overlooked in the
literature.}}

\begin{theorem} \cite{KZ} 
The billiard in a typical polygon is
topologically transitive.
\end{theorem}

\begin{proof}
Let $\D \subset \R^2$ be the disc.
Identify the phase space of the billiard flow in each $n$-gon
with $\D \times \S$ and assume that this identification depends
continuously on the polygon.
Let $B_i$ be a countable basis for the topology of  $\D \times \S$.

Fix an nonempty open set $U_k \subset \D \times \S$ (the choice of this set
will be made precise below) and let $K_k \subset K$
be the set of $n$-gons $P$ such that 
there exists a billiard trajectory starting in $U_k$
that visits all (the images of) the sets $B_1,\dots,B_k$ in the phase space
of the billiard flow in $P$.  Each $K_k$ is open, and their intersection
is a $G_{\delta}$ set.

To see denseness let $Y_q$ be the set of rational $n$-gons with angles
$\pi p_i/q_i$, with $p_i,q_i$ co-prime and the least common multiple 
of the $q_i$ at least $q$.
For every $P \in Y_q$ each invariant surface $M_\theta$ is $1/q$ dense
in the phase space. Therefore for every $k$ there exists $q$ such that
for every $P \in Y_q$ the surface $M_{\theta}$ (for all $\theta$) 
intersects $U_k$ and all the (images of) 
the sets $B_1,\dots,B_k$ in the phase space of the billiard flow in $P$.
Since $M_\theta$ has a dense trajectory for all but countably many
$\theta$ we have $Y_q \subset X_k$. Thus since $Y_q$ is dense in 
the space of $n$-gons, the set $K_k$ is as well.  It follows that $\cap K_k$ is
dense.

Let $P$ be a polygon in $\cap K_k$.  The choice of $U_1$  is
arbitrary. Let $\tilde{U}_1 \subset U_1$ be a nonempty
compact neighborhood. 
Suppose we have inductively
chosen a nonempty compact neighborhood 
$\tilde{U}_{k-1} \subset U_{k-1}$. 
Let $U_k \subset \tilde{U}_{k-1}$ be open.
Since $P \in K_k$
we can find a billiard trajectory starting in $U_k$ that visits each
$B_1,\dots,B_k$. By continuity there is an compact neighborhood $U_k \subset
\tilde{U}_{k-1}$ such that each trajectory in $U_k$ visits $B_1,\dots,B_k$ within
a bounded time $T_k$.  Thus any $x \in \cap U_k = \cap \tilde{U}_k$ has a dense trajectory,
possibly singular.\footnote{This detail is also often overlooked in
the literature.}

To see that there are non-singular trajectories it suffices to 
produce an uncountable 
collection of such dense trajectories.  Since there are only countably
many saddle connections, i.e.~orbits segments starting and ending
at a vertex, most of the dense trajectories are not saddle connections.
If such a trajectory is singular, i.e.~the $n$th iterate arrives at a
vertex, then the (forward) orbit starting at time $n+1$ is dense and non-singular.

Thus we modify the construction to produce a Cantor set 
(an uncountable set) of points with dense trajectories. To do this 
choose $\tilde{U}_{k}$ be the disjoint union of $2^{k}$ nonempty
compact neighborhoods such that each neighborhood contains exactly
two of the neighborhoods of  $\tilde{U}_{k+1}$ and such that there
is a billiard trajectory in each of the neighborhoods which visits the
sets $B_1,\dots,B_k$.  
\end{proof}

In 1986, Kerckhoff, Masur and Smillie, using Teichmüller theory,
proved that the directional billiard flow is 
uniquely ergodic for almost every direction \cite{KMS}. 
The surfaces $M_{\theta}$ are not only
$1/q$ dense in the phase space, but also approximately well distributed.
Combining these two facts and  no other new ideas one can conclude:

\begin{theorem}
\cite{KMS,K1,PS}
The billiard in a typical polygon is 
ergodic.
\end{theorem}

Let $\phi:\N \to \R^+$ such that
$\lim_{q \to \infty} \phi(q) = 0$. Let $P$ be a $n$-gon with angles
$\alpha_1,\dots,\alpha_n$.  We say that $P$ admits approximation by
rational polygons  at rate $\phi(q)$ if for every $q_0 >0$ there is $q > q_0$
and positive integers $p_1,\dots,p_n$ each co-prime with $q$ such that 
$|\alpha_i - \pi p_i/q| < \phi(q)$ for all $i$.

Vorobets has proven the following constructive ergodicity result:
\begin{theorem} \cite{V}
Let $P$ be a polygon that admits approximation by rational
polygons at the rate $$\phi(q) = \left (2^{2^{2^{2^q}}}\right )^{-1},$$
then the billiard flow in $P$ is ergodic.
\end{theorem}

I proved the following theorem:
\begin{theorem}\label{total}\cite{Tr2004}
The billiard map in a  typical polygon is totally ergodic,
i.e.~$T^n$ is ergodic for all $n \in \N$.
\end{theorem}

Both the constructive ergodicity and the total ergodicity 
theorems are proven by using
an alternative version of the proof of the result of Kerckhoff, Masur and Smillie
based on a combination of ideas of Masur and Boshernitzan \cite{M,Bo2}.

In 1999 I applied the approximation
method to infinite polygons \cite{Tr1999}.  My result was improved by
Degli Esposti, Del Magno and Lenci
who showed:

\begin{theorem}\cite{DDL} 
Suppose $p_n$ is a monotonically decreasing
sequence of positive numbers satisfying $\sum p_n  = 1$.
Let $$P = \cup_{n \ge 0} [n,n+1] \times [0,p_n].$$ 
Then the billiard in a typical $P$ is ergodic for a.e.~direction.
\end{theorem}

Here the topology on the space of such $P$ is given by the metric
$$d(P,Q) := \sum |p_n - q_n| = Area(P \triangle Q).$$

To prove this result  approximate $P$ by finite polygons
from the class, i.e.~$q_n \equiv 0,$ for $n \ge N$.

\section{Old approximation results for convex 
smooth tables.}

Let $C$ be a strictly convex billiard table. A compact convex set $K$ 
is a {\em caustic}  of $C$ if the boundary of $C$ is obtained by wrapping 
a string around $K$, pulling it tight at a point and moving the point 
around $K$ while keeping the string tight.
If a billiard orbit is tangent to a caustic once, then it
is tangent to it in-between every pair of bounces.

In 1979 Lazutkin showed: 
 
\begin{theorem}\cite{L}
If $C$ is a sufficiently smooth strictly convex billiard table 
then the table contains ``many'' caustics.
\end{theorem}

Many means the union of the caustics has positive area.

\begin{corollary}
The billiard in a sufficiently smooth strictly convex
table is not ergodic, not even topologically transitive.
\end{corollary}

The proof of Lazutkin's result is based on KAM theory.  
Lazutkin's original proof need the table
to be $C^{553}$.  
Rüssmann's version of KAM can be used to replace 553 by
8, and finally R.~Douady's version of KAM shows that 7 is sufficient \cite{Ru,D}.

In 1990 Gruber noticed that we can apply the approximation ideas
of the previous section to prove contrasting results for 
the low smoothness case.

\begin{theorem}\cite{G}
The $C^0$-typical convex billiard table has the following properties:
\begin{enumerate}
\item it contains no caustics,
\item it is strictly convex,
\item it is of class $C^1$,
\item it is topologically transitive.
\end{enumerate}
\end{theorem}

To prove the topological transitivity we approximate a convex table
by finite polygons with increasing number of sides.

In fact, using approximation teachings and the results of 
Kerckhoff, Masur and Smillie
one can easily conclude that

\begin{theorem}
$C^1$-typical
$C^1$-convex billiard table is strictly convex and ergodic.
\end{theorem}

\section{New results on smooth tables.}

With A.~Stepin we have have improved Gruber's result.
To describe our result we begin by a new characterization of weak
mixing.

\subsection{Weak mixing.}
Let $(X,\beta,\mu)$ be a probability space and $T: X \to X$ a measure
preserving transformation.  $T$ is {\em weak mixing} iff
$$\forall A,B \in \beta: \ \lim_{n \to \infty} \frac1n \sum_{k=0}^{n-1}|\mu(T^{-k}A \cap B) - \mu(A)\mu(B)| = 0$$
iff $\forall A,B \in \beta$: there exists $J= J(A,B)  \subset \N$ of full 
density (i.e.~$(\#n\in J: n \le N)/N \to 1$) such that
$$\quad \lim_{n \in J \to \infty} \mu(T^{-k}A \cap B) =\mu(A)\mu(B).$$

A set $J \subset \N$ is called {\em full} if $J$ contains arbitrarily long runs,
i.e.~$\forall N, \exists k$ s.t.~$\{k,k+1,\dots,k+N\} \subset J$.

\begin{theorem}\cite{ST}
$T$ is weak-mixing iff $\forall A \in \beta$
there exists $J$ full such that
$$\lim_{n \in J \to \infty} \mu(T^{-k}A \cap A) =\mu(A)^2$$
\end{theorem}
\begin{proof} 
First we claim that $T$ is ergodic.  Suppose $T^{-1}(A) = A$ mod $0$.  
By invariance we have $\mu(T^{-n}A \cap A) = \mu(A)$, but by 
mixing along full sequences this quantity converges to $\mu(A)^2$ for
$n \in J$.  Thus $\mu(A) = 0$ or $\mu(A) =1$.

Now we are ready to prove that $T$ is weak-mixing.
Suppose not, then there exists a Kronecker factor  $(\hat{X},\hat{\beta},\hat{\mu},\hat{T})$, i.e.~a rotation
$\hat{T}$ on a compact Abelian group $\hat{X}$ such that the following
diagram commutes:
$$X\mrfn{T}X$$
$$\dfn{\pi} \quad \ \  \dfn{\pi}$$
$$\hat{X}\mrfn{\hat{T}}\hat{X}$$

Fix $\hat{A} \in \hat{\beta}$ such that $0 < \hat{\mu}(\hat{A})< 1$ and
let $A := \pi^{-1}(\hat{A})$.  

Then since $\hat{T}$ is a rotation, there is a quasi-periodic sequence $n_i$
such that $\hat{\mu}(\hat{T}^{-n_i}\hat{A} \cap \hat{A}) 
\approx\hat{\mu}(\hat{A})$ and thus
$$\mu(T^{-n_i}A \cap A) = \hat{\mu}(\hat{T}^{-n_i}\hat{A} \cap \hat{A}) 
\approx\hat{\mu}(\hat{A}) = \mu(A)$$
which contradicts fullness.
\end{proof}

\subsection{Hyperbolic billiards.}

Instead of approximating a smooth table by polygons we approximate
by certain hyperbolic tables. Namely
consider a polygon $P$.  Replace each vertex with a (small) circular
arc such that 1) the resulting table is $C^1$ and 2) the focusing circles
lie inside the table. We call such tables B-tables and if the original
polygon is convex then CT-tables.

\begin{theorem}
(Bunimovich \cite{B}) For any  B-tables 
each ergodic component is open (mod 0) and the billiard on each
ergodic component is mixing, K-mixing, Bernoulli.
\end{theorem}

\begin{theorem}
(Chernov, Troubetzkoy \cite{CT}) Any CT-table is ergodic, hence
mixing,  K-mixing, Bernoulli.
\end{theorem}

\subsection{Back to  $C^1$ billiards.}
\begin{theorem}\cite{ST}
The billiard flow/map in a $C^1$-typical $C^1$ billiard table is weak-mixing.
\end{theorem}

{\bf Sketch of the proof:} 
If the table is convex we approximate by billiards from the class of
CT-tables introduced
above.  Fix a CT-table and a
finite collection of sets.  For all tables sufficiently close to the fixed 
table we can find an arbitrarily long run of times such that 
these sets approximately mix along this run.

If the table is not convex then we do not have a nice class of ergodic tables
like the CT-class. 
To prove the theorem we approximate the table
by polygons to obtain generic ergodicity.  Then we approximate
by B-tables to get the weak-mixing.\hfill\qed
\subsection{Convex tables.}
We can slightly improve the smoothness is the table is convex, or piecewise
concave-convex. Let PCC be the class of $C^1$ piecewise convex-concave
billiard domains.
Because of the piecewise monotonicity of the derivative, 
locally the second derivative exists almost everywhere.
Consider the Hausdorff distance $\rho_H$ between the second derivative
of the boundaries of such tables. 
Let $B^2$ be the closure of the class
of B-tables in the topology induced by the metric
$\rho = \rho_H + \rho_{C^1}$.

\begin{theorem}\cite{ST}
Weak-mixing is typical in $B^2$
\end{theorem}
We remark that 
$PCC \subset B^2 \subset C^1$ and
that every table in $PCC$ is $C^2$ except at a finite number of points.
Furthermore the set of strictly convex tables is residual in $B_{convex}^2$, thus
the generic table in $B_{convex}^2$ is strictly convex and weak-mixing.

\section{Back to polygonal billiards.}
Recently Avila and Forni have shown that a.e.~interval
exchange not of rotation type is weak-mixing \cite{AF}. 
Furthermore, for almost every
translation surface they can show that the directional flow is weak-mixing.
A priori, these results do not hold for polygonal billiards. Thus the
results obtained by Stepin and myself stated here are conditional.

Let $WMix_c$ be the set of rational polygons for which at least $c\%$ of the ergodic components are weak-mixing.  Note that the square and the equilateral
triangle are not in $WMix_c$ for any $c > 0$. The result of Avila-Forni gives
hope that all but
finitely many rational polygons are in $WMix_1$.

\begin{theorem}\cite{ST}
1) If there exists $c >0$ such that $WMix_c$ is 
dense in the set of polygons then weak-mixing is generic for polygonal 
billiards.\\
2) Either there exists $m$ such that $WMix_{\frac1m}$ is somewhere dense in
the set of polygons (and hence weak-mixing is on 2nd category for polygonal
billiards) or $\cup_m WMix_{\frac1m}$ is nowhere dense.
\end{theorem}

\section{Angular recurrence and periodic billiard orbits.}

The billiard orbit of any point which begins perpendicular to a side 
of a polygon and at a later instance hits some side perpendicularly 
retraces its path infinitely often in both senses between the two 
perpendicular collisions and thus is periodic.

\subsection{History and notation.}
In 1991 Ruijgrok conjectured  based on numerical evidence that for any 
irrational triangle almost every orbit which starts perpendicular to a 
side is twice perpendicular and thus periodic\cite{R}.
Here almost every is with respect to the length measure on the side.
In 1992 Boshernitzan \cite{Bo} and independently Galperin, Stepin and
Vorobets
\cite{GSV} proved this conjecture for all rational polygons.
In this case the set of non-periodic perpendicular orbits is
a finite set.
In 1995 Cipra, Kolan and Hansen proved this conjecture for all
right triangles \cite{CHK}.
Apparently these four articles are completely independent of each other.

Applying Theorem \ref{total}, I obtained the first
approximation result in this direction.

\begin{theorem}\cite{Tr2004}
The conjecture is true for dense $G_{\delta}$ of $n$--gons
with $n-2$ angles fixed rational multiples of $\pi$.
\end{theorem}

Billiards are not affected by similarity, thus it is
convenient to fix two corners, say at the origin and the
point $(1,0)$. For triangles if we fix one angle to be a rational
multiple of $\pi$ then the set of triangles is an interval
parametrized by one of the other angles. This theorem says
that for a dense $G_{\delta}$ subset of (any compact subset of)
this interval the conjecture is true.

To discuss our results on billiards in parallelograms we must
introduce some more notation.
For the moment we assume that $P$ is convex.
We will describe
the billiard map $T$ as a transformation of the set $X$ 
of rays which intersect $P$.
Let  $\t$ be the angle between the billiard trajectory
and the positive $x$--axis. 
Consider the perpendicular cross section 
$\Xt$ to
the set of rays whose angle is $\t$.  The set $\Xt$ is
simply an interval.
Let $w$ be the arc-length on $\Xt$. 

For a non-convex polygon we must differentiate the portion
of rays which enter and
leave the polygon several times. The above construction can be done
locally, yielding the set $\Xt$ which consists
of a finite union of intervals.
Let $w$ be the unnormalized length measure on $\Xt$.

For every direction $\t$ let 
$\Bt$ be the set of $x \in \Xt$ whose forward orbit never returns parallel
to $x$.
A polygon is called {\em angularly recurrent} if for every direction $\t$ we have
$w(\Bt) = 0$.
Angularly recurrent polygons satisfy the conclusion of Ruijgrok's
conjecture, almost every perpendicular orbit is periodic!
There is no reason to think that a polygon chosen at random should
be angularly recurrent, this property is much stronger than Poincaré
recurrence which only implies that almost every point comes back close
to itself (angles close not equal!).

A polygon $P$ is called a {\em generalized parallelogram} if all sides of $P$
are parallel to two fixed vectors.  Note that a right triangle unfolds
to a rhombus and thus results for generalized parallelograms hold
also for right triangles.

In 1996 Gutkin and I generalizing the result of Cipra, Kolan and
Hansen proved that all right triangles and all generalized 
parallelograms  are angularly recurrent \cite{GT}.

\subsection{Dimensions.}
Let $Y \subset \bf{R^n}$.  Let $N(\epsilon)$ denote 
the minimal number of $\epsilon$ balls needed to cover $Y$.  
The {\em lower box dimension} of $Y$, 
denoted by $\dim_{LB}Y$ is given by
$$
\liminf_{\epsilon \to 0} \frac{\log N(\epsilon)}{\log 1/\epsilon} \ .
$$
The {\em upper box dimension} $\dim_{UB}$ is defined similarly, replacing
the $\liminf$ by $\limsup$.  
If $\dim_{UB}Y$ and $\dim_{LB}Y$ both exist and are equal, we define 
the box dimension of $Y$ to be this value, and write $\dim_{B}Y$ = 
$\dim_{UB}Y$ = $\dim_{LB}Y$.  

Let $s \in [0,\infty]$.  The $s$-{\em dimensional 
Hausdorff measure} ${\mathcal H}^s (Y)$
of a subset $Y$ of $\bf{R^n}$ is defined by the following limit 
of covering sums:
\begin{eqnarray*}
{\mathcal H}^s (Y) & = & \lim_{\epsilon \to 0} ( \inf
\{ \sum^\infty_{i=1} \left(\hbox{ diam } U_i\right)^s :\\
&& Y \subset \bigcup^\infty_{i=1} U_i \hbox{ and } \sup_{i} 
\hbox{ diam } U_i \leq \epsilon  \}  ).
\end{eqnarray*}

It is easy to see that there exists a unique $s_0 = s_0 (Y)$ such that
$$
{\mathcal H}^s (Y) = \infty \mbox{ for } s < s_0 \mbox{ and }
  0 \mbox{ for } s > s_0\ .$$

The number $s_0$ 
is called the {\em Hausdorff dimension} of $Y$ and is
denoted by $\dim_H Y$.
Standard arguments give that 
$$
\dim_H Y \leq \dim_{LB} Y\leq \dim_{UB} Y
$$

\subsection{Approximation questions.}

With Jörg Schmeling we showed that for 
all generalized parallelograms $\dim_{LB}(F_\t) \le 1/2$ for all $\t$ \cite{ScTr}.

The following theorem shows that we
can do better that 1/2 for well approximable directions in a given
parallelogram or for all directions in a  parallelogram with well
approximable angle.

\begin{theorem}\label{thm2}\cite{ScTr}
Fix a generalized parallelogram with irrational angle $\a$.
For any $\t$ such that $\| (\t+ p\a)/\pi\| < p^{-\mu}$
has infinitely many  solutions $p \in \N$. Then 
$$\dim_{LB} (F_\t \cap \Ut) \le \frac{1}{\mu+1}.$$
Similarly $\dim_{LB} (F_\t \cap \Dt) \le \frac{1}{\mu+1}$
provided the angle $\t - \a/2$ is $\mu$-well approximable.
\end{theorem}
 
Two questions naturally arise,
{\em for a fixed generalized parallelogram how many directions have small
$F_\t$} and
{\em for a fixed direction $\t_0$ how many parallelograms have a small $F_{\t_0}$?}

Our answer to these question required new and recent 
results in diophantine approximation.

\subsection{New number theoretic results.}

For $\mu \ge 1$ and $\w \in \S$ let 
\begin{eqnarray*}\Am  :=  \{ t \in \S: \quad \| t+ p\w\| < p^{-\mu}
 \hbox{ for infinitely many } p \in \N\}.
\end{eqnarray*}

Minkowski's classical theorem states that if $\mu = 1$ then 
$\Am$ consists of all $t \in \S$ such
that $t$ is not in the orbit of $\w$.

Using different techniques Schmeling and I and independently Y.\ Bugeaud
have shown: 
\begin{theorem}\label{thm3}\cite{ScTr},\cite{Bug}
For any $\mu > 1$ we have
$$\dim_H (\Am) = 1/\mu.$$
\end{theorem}

For $\mu \ge 1$ and $t \in \S$ let 
\begin{eqnarray*}
\Bm := \{\w \in \S: \|t+p\w\| < p^{-\mu}
\hbox{ for infinitely many } p \in \N\}.
\end{eqnarray*}

The other result we need to apply is a recent one of Levesly:

\begin{theorem} \label{thmLe}(Levesly \cite{Le})
For any $\mu > 1$ we have
$$\dim_H(\Bm) = \frac{2}{1+\mu}.$$ 
\end{theorem}
Remark: this is a classical theorem of Jarnik for $t=0$.

Besides the applications to billiards these results can be 
used to analyze the rotation driven
dynamical variant of
the classical Dvoretzky covering of the circle \cite{FS}.
More recently we have studied dynamical diophantine approximation
for the circle doubling map \cite{FST}.

\subsection{Application to billiards.}

We reply to the first question, for 
a fixed generalized parallelogram how many directions have small
$F_\t$?  Let
$$\Cs := \{\t: \dim_{LB}(\Bt \cap \Ut) \le s\}.$$
\begin{theorem}\cite{ScTr}
For all $s \in [0,1/2]$ 
$$\dim_H(\Cs) \ge \frac{s}{1-s}.$$
The set $\Cs$ is residual with box dimension 1.
\end{theorem}
\begin{proof}
This follows immediately from Theorems \ref{thm2} and \ref{thm3}.
\end{proof}
 
We turn to the second question, for a fixed direction $\t_0$ how many parallelograms have a small $F_{\t_0}$?
Let 
\begin{eqnarray*}
\Ds & := & \{\a \in \S: \dim_{LB}(F_{\t_0} \cap U_{\t_0} )<s \hbox{ for all }\\
&&\hbox{ generalized  parallelograms with angle } \a  \}.
\end{eqnarray*}
\begin{theorem}\cite{ScTr}
For $\t_0$ fixed and $s \in [0,1/2]$ we have 
$$\dim_H\Ds \ge 2s.$$
The set $\Ds$ is residual and has box dimension 1.
\end{theorem}
\begin{proof}
This follows immediately from Theorems \ref{thm2} and \ref{thmLe}.
\end{proof}

\end{document}